\documentclass[12pt]{article}
\usepackage{graphicx}
\usepackage{amsmath,amsthm,amssymb,enumerate}%, esint}
\usepackage{euscript,mathrsfs}
\usepackage{color}
\usepackage{dsfont}
\usepackage{multirow}
\usepackage{multicol}
\usepackage{url}
\usepackage{esint}
\usepackage[left=2cm,right=2cm,top=3.5cm,bottom=3.5cm]{geometry}
\usepackage{color}
\usepackage[framemethod=tikz]{mdframed}
\allowdisplaybreaks

\usepackage{soul}

\catcode`\@=11 \@addtoreset{equation}{section}

\catcode`\@=12

\usepackage{subcaption}

\newtheorem{Theorem}{Theorem}[section]
\newtheorem{Proposition}[Theorem]{Proposition}
\newtheorem{Lemma}[Theorem]{Lemma}
\newtheorem{Corollary}[Theorem]{Corollary}

\theoremstyle{definition}
\newtheorem{Definition}[Theorem]{Definition}

\newtheorem{Remark}[Theorem]{Remark}

\newcommand{\bTheorem}[1]{
	\begin{Theorem} \label{T#1} }
	\newcommand{\eT}{\end{Theorem}}

\newcommand{\bProposition}[1]{
	\begin{Proposition} \label{P#1}}
	\newcommand{\eP}{\end{Proposition}}

\newcommand{\bLemma}[1]{
	\begin{Lemma} \label{L#1} }
	\newcommand{\eL}{\end{Lemma}}

\newcommand{\bCorollary}[1]{
	\begin{Corollary} \label{C#1} }
	\newcommand{\eC}{\end{Corollary}}

\newcommand{\bRemark}[1]{
	\begin{Remark} \label{R#1} }
	\newcommand{\eR}{\end{Remark}}

\newcommand{\bDefinition}[1]{
	\begin{Definition} \label{D#1} }
	\newcommand{\eD}{\end{Definition}}

\newcommand{\Del}{\Delta_x}
\newcommand{\tvu}{\widetilde{\vu}}

\newcommand{\avintQ}[1]{\fint_{Q} #1 \dx}

\newcommand{\prst}{\mathbb{P}}
\newcommand{\tTh}{\widetilde{\Theta}}

\newcommand{\tvT}{\tilde{\Theta}}

\newcommand{\bFormula}[1]{
	\begin{equation} \label{#1}}
	\newcommand{\eF}{\end{equation}}

\newcommand{\vuh}{\vu_h}

\newcommand{\Ov}[1]{\overline{#1}}

\newcommand{\aleq}{\stackrel{<}{\sim}}

\newcommand{\vt}{\vartheta}
\newcommand{\vu}{\vc{u}}

\newcommand{\vc}[1]{{\bf #1}}

\newcommand{\Div}{{\rm div}_x}
\newcommand{\Grad}{\nabla_x}

\newcommand{\dx}{\,{\rm d} {x}}

\newcommand{\dt}{\,{\rm d} t }

\newcommand{\intQ}[1]{\int_{{Q}} #1 \ \dx}

\newcommand{\D}{{\rm d}}

\newcommand{\bS}{\mathbb{S}}

\newcommand{\expe}[1]{ \mathbb{E} \left[ #1 \right] }

\newcommand{\br}{ \nonumber \\ }

\def\softd{{\leavevmode\setbox1=\hbox{d}%
		\hbox to 1.05\wd1{d\kern-0.4ex{\char039}\hss}}}
\definecolor{Cgrey}{rgb}{0.85,0.85,0.85}
\definecolor{Cblue}{rgb}{0.50,0.85,0.85}
\definecolor{Cred}{rgb}{1,0,0}
\definecolor{fancy}{rgb}{0.10,0.85,0.10}
\definecolor{amaranth}{rgb}{0.9, 0.17, 0.31}

\newcommand\Cbox[2]{%
	\newbox\contentbox%
	\newbox\bkgdbox%
	\setbox\contentbox\hbox to \hsize{%
		\vtop{
			\kern\columnsep
			\hbox to \hsize{%
				\kern\columnsep%
				\advance\hsize by -2\columnsep%
				\setlength{\textwidth}{\hsize}%
				\vbox{
					\parskip=\baselineskip
					\parindent=0bp
					#2
				}%
				\kern\columnsep%
			}%
			\kern\columnsep%
		}%
	}%
	\setbox\bkgdbox\vbox{
		\color{#1}
		\hrule width  \wd\contentbox %
		height \ht\contentbox %
		depth  \dp\contentbox
		\color{black}
	}%
	\wd\bkgdbox=0bp%
	\vbox{\hbox to \hsize{\box\bkgdbox\box\contentbox}}%
	\vskip\baselineskip%
}

\mdfdefinestyle{MyFrame}{%
	linecolor=black,
	outerlinewidth=1pt,
	roundcorner=5pt,
	innertopmargin=\baselineskip,
	innerbottommargin=\baselineskip,
	innerrightmargin=10pt,
	innerleftmargin=10pt,
	backgroundcolor=white!20!white}

\newcommand{\red}{\textcolor{red}}

\newcommand{\Gradh}{\nabla_h}

\newcommand{\vtB}{\vt_{B}}

\newcommand{\difuh}{\bS_h:\Gradh \vuh }

\newcommand{\R}{\mathbb{R}}

\newcommand{\bfphi}{\boldsymbol{\varphi}}

\allowdisplaybreaks

\begin{document}

%%%%%%%%%%%%%%%%%%%%%%%%%%%%%%%%

\title{\bf A continuous data assimilation method for 
a variant of Oberbeck--Boussinesq system with randomly 	
perturbed data}

\author{Eduard Feireisl$^{1}$
	\thanks{
		The work of E.F.\ was partially supported by the
		Czech Sciences Foundation (GA\v CR), Grant Agreement
		24--11034S. The Institute of Mathematics of the Academy of Sciences of
		the Czech Republic is supported by RVO:67985840.
		E.F.\ is a member of the Ne\v cas Center for Mathematical Modelling.} 
	\and  M\u ad\u alina Petcu$^{2,3,4}$  }

	%\date{\today}
	
	\maketitle

	\medskip
	
	\centerline{${}^1$Institute of Mathematics of the Academy of Sciences of the Czech Republic}
	\centerline{\v Zitn\' a 25, CZ-115 67 Praha 1, Czech Republic}
	
	\bigskip
	
	\centerline{${}^2$Laboratoire de Math\' ematiques et Applications, UMR CNRS 7348 - SP2MI}
	
	\centerline{Universit\' e de Poitiers, Boulevard Marie et Pierre Curie - T\' el\' eport 2}
	
	\centerline{86962 Chasseneuil, Futuroscope Cedex,
		France}
	
	\centerline{${}^3$The Institute of Mathematics of the Romanian Academy, Bucharest, Romania}

\date{}

\maketitle

\medskip

\begin{abstract}
	
We show convergence of a continuous data assimilation method for the 
Oberbeck-Boussinesq system in the dimension $d=2,3$. Our working hypothesis is boundedness of the reference solution, while the synchronized solution satisfies the equations in a weak sense. The main tool  is the relative energy inequality for stochastic problems.

\end{abstract}

%\bigskip

{\small
	
	\noindent
%	{\bf 2020 Mathematics Subject Classification:}{ 
%		(primary); 
%		(secondary) }
	
	\medbreak

	%\tableofcontents
	
}

\section{Introduction}
\label{i}

We consider a variant of the Oberbeck--Boussinesq (OB) system identified 
in \cite{BelFeiOsch} as a singular limit of a mildly stratified 
Navier--Stokes--Fourier system. The time evolution of the 
velocity $\vu = \vu(t,x)$ and the temperature fluctuations $\theta = \theta(t,x)$ is described by the following Oberbeck--Boussinesq (OB) system of equations:
\begin{align}
\Div \vu &= 0 , \br 
\partial_t \vu + \Div (\vu \otimes \vu) + \Grad \Pi &= 
\mu \Del \vu - \theta \Grad G,\ \mu > 0,  \br	
\partial_t \theta + \Div (\vu \theta) + a \Div (G \vu) &= \kappa \Del \theta, \ a \in \R, \ \kappa > 0,
\label{i1}	
\end{align}
considered in a bounded domain $Q \subset \R^d$, and supplemented with the boundary conditions 
\begin{align}
	\vu|_{\partial Q} &= 0, \br 
\theta|_{\partial Q} &= \theta_B - \alpha \avintQ{\theta },\ \alpha \in [0,1). 	
\label{i2}
\end{align}	
$G = G(x)$ represents a gravitational potential, $\avintQ{ }= \frac{1}{|Q|} \int_Q$ and $\theta_B=\theta_B(x)$ is the boundary distribution of the original temperature deviation. The value of the parameter $\alpha$ depends on the reference state of the system. The non-local boundary condition \eqref{i2} can be seen as a proper singular limit of 
the Dirichlet boundary conditions imposed on the original system before scaling, see \cite{BelFeiOsch}, \cite{Fa-Fe_2024}. 

Given its physical 
background, the OB system represents a simple mathematical model with potential applications in meteorology. Motivated by the seminal work 
of Azouani, Olson, and Titi \cite{AzOlTi}, we consider the problem of continuous data assimilation for system \eqref{i1}, \eqref{i2}. Similarly 
to Bessaih, Olson, and Titi \cite{BeOlTi}, we implement appropriate nudging operators producing a possible stochastic error. Accordingly, the 
resulting synchronized problem consists of a system of PDE's with stochastic forcing. 

Our strategy is based on the relative energy inequality 
developed in the context of stochastic PDE's in \cite{BrFeHo2015A},  see also \cite{FeiPet}:
\begin{itemize}
\item	
Under a mild assumption of boundedness, the solutions of the observed OB system are eligible as test functions in the relative energy inequality  
for the associated synchronized system.

\item Performing standard stochastic analysis, we obtain qualitative estimates on the distance between the observed and synchronized solutions 
in terms of expected values.  

\end{itemize}
\noindent
We point out that the method works in the physically relevant case $d=3$, where the well posedness for the observed system remains an outstanding open problem.

The paper is organized as follows. Section \ref{P} contains the necessary preliminary material as well as 
a suitable reformulation of the observed OB system. In Section \ref{r}, we formulate the data assimilation problem and introduce the 
associated synchronized system. The main result is stated in Section \ref{M}. In Section \ref{e}, we introduce the relative energy inequality, and, 
finally, complete the proof of the main result in Section \ref{a}.

\section{Preliminaries, problem formulation}
\label{P}

We suppose $Q \subset \R^d$, $d=2,3$ is a bounded domain of class $C^{2,\red{\beta}}$, $\beta > 0$. The regularity 
of the boundary can be slightly relaxed but still needed as the parabolic regularity for the Stokes problem will be applied. 

Without lost of generality, we may suppose 
\begin{equation} \label{i3}
\intQ{ G } = 0,
\end{equation}
and also 
\begin{equation} \label{i4}
\Del G = 0, 
\end{equation}
which seems a reasonable assumption for a gravitational potential.

Finally, we denote $\vtB$ the harmonic extension of 
the boundary data $\theta_B$ inside $Q$: 
\begin{equation} \label{i5}
\Del \vtB = 0 \ \mbox{in}\ Q,\ \vtB|_{\partial Q} = \theta_B.	
\end{equation}

\subsection{Reformulation of the problem}

For future analysis, it is convenient to reformulate the OB system to a problem with \emph{homogeneous} Dirichlet boundary conditions.
Setting 
\[
\Theta = \theta + \alpha \avintQ{ \theta } - \vtB
\]
we may rewrite system \eqref{i1}, \eqref{i2} in the form
\begin{align}
	\Div \vu &= 0 , \br 
	\partial_t \vu + \Div (\vu \otimes \vu) + \Grad \Pi &= 
	\mu \Del \vu - \Theta \Grad G - \vtB \Grad G,  \br	
	\partial_t M(\Theta) + \Div (\vu \Theta)  &= \kappa \Del \Theta - \Div (\vu (\vtB + aG)) ,
	\label{i6}	
\end{align}
with the boundary conditions
\begin{align}
	\vu|_{\partial Q} = 0, \
	\Theta|_{\partial Q} = 0.
	\label{i7}
\end{align}	
The operator $M$ is defined as 
\begin{equation} \label{i8}
M(\Theta) = \Theta - \frac{\alpha}{\alpha + 1}\avintQ{ \Theta}.
\end{equation}	

It is straightforward to check that 
\[
M : L^2(Q) \to L^2(Q) 
\]
is a bounded, self--adjoint, positively definite operator, 
with a bounded inverse 
\begin{equation} \label{i9}
M^{-1} (h) = h + \alpha \avintQ{ h },\ h \in L^2(Q).
\end{equation} 

We call system \eqref{i6}--\eqref{i8} \emph{modified Oberbeck--Boussinesq (MOB) system}. 

\subsection{Mathematics of the MOB system}
\label{mat}

The existence of global--in--time weak solutions for 
the MOB system was established in \cite{AbbFei22}. In particular, 
there exists a global in time weak solution in the class 
\begin{align} 
\vu &\in L^\infty (0, T; L^2(Q; \R^d)) \cap L^2(0,T; W^{1,2}_0(Q; \R^d)), \br 
\Theta &\in L^\infty((0,T) \times Q) \cap L^2(0,T; W^{1,2}_0(Q))	
\label{i10}
\end{align}
for any initial data 
\begin{align}
\vu(0, \cdot) = \vu_0 &\in L^2(Q; \R^3),\ \Div \vu_0 = 0,\ \vu_0 \cdot \vc{n}|_{\partial Q} = 0, \br 
\Theta(0, \cdot) = \Theta_0 &\in C(\Ov{Q})
\label{i11}	
\end{align}
provided 
\begin{itemize}
	\item $Q \subset \R^d$ is a bounded Lipschitz domain; 
	\item $G \in W^{1,\infty}(Q)$; 
	\item $\vtB \in C(\Ov{Q}) \cap W^{1,2}(Q)$,  
\end{itemize}
see \cite[Theorem 2.3]{AbbFei22}. 

If $d=2$, the weak solution is unique and regular as long as the domain 
as well as the data enjoy additional regularity. Specifically, if  
\begin{itemize}
	\item $Q \subset \R^2$ is a bounded domain of class $C^2$; 
	\item $\vtB \in C^2(\Ov{Q})$;  
	\item $\vu_0 \in W^{2, \infty}(Q; \R^2)$ satisfies the compatibility condition $\vu_0|_{\partial Q} = 0$; 
	\item $\Theta_0 \in W^{2,\infty}(Q)$ satisfies the compatibility condition 
	$\Theta_0|_{\partial Q} = 0$;
\end{itemize}
then the MOB problem admits a unique global in time strong solution in the class 
\begin{align}
\vu &\in L^q(0,T; W^{2,q}(Q; \R^2)),
\partial_t \vu \in L^q(0,T; L^{q}(Q; \R^2)), \br 
\Theta &\in L^q(0,T; W^{2,q}(Q; \R^2)),
\partial_t \Theta \in L^q(0,T; L^{q}(Q; \R^2))
\label{i12}
\end{align}
for any finite $q \in [1, \infty)$, see \cite[Theorem 3.1]{AbbFei22}.

Moreover, as shown in \cite[Section 4]{FeRoSc2024}, the MOB system 
is dissipative in the sense of Levinson and admits a trajectory attractor $\mathcal{A}$ that consists of all \emph{entire solutions}: 
\begin{align} 
\mathcal{A} &= \left\{ (\vu, \Theta) \ \Big|\ 
(\vu, \Theta) \ \mbox{weak solution of MOB system for}\ t \in (- \infty, \infty),\ 
\right. \br 
&\quad \quad \sup_{t \in \R} \left( \| \vu(t, \cdot) \|_{L^2(Q; \R^2)} + 
\| \Theta (t, \cdot) \|_{L^\infty(Q)} \right) \leq \mathcal{F} \Big\},
\label{i13}
\end{align}
where $\mathcal{F}$ depends only on $\alpha$, $\max_{\partial Q}|\vt_B|$, the diffusion coefficients and the Poincar\'e constant, but is independent of the initial conditions.

\section{Observed and synchronized solutions}
\label{r}

We introduce the concept of reference and synchronized solution
to the MOB problem.

\subsection{Observed solution}

The \emph{observed solution} $(\vu, \Theta)$ is a solution of the 
MOB system defined on a finite time--interval,  
\[
t \in (T^-, T^+),\ - \infty < T^- < 0 < T^+ < \infty
\]
satisfying
\begin{equation} \label{r1}
\sup_{t \in (T^-, T^+)} \left( \| \vu (t, \cdot) \|_{L^\infty(Q; \R^d) } + 
\| \Theta (t, \cdot) \|_{L^\infty(Q) }\right) \leq \mathcal{E}.  	
\end{equation}

In view of the existence results recorder in Section \ref{mat}, 
any solution belonging to $\mathcal{A}$ satisfies \ref{r1} if $d = 2$. 
Moreover, in view of the regularizing effect of the Laplace/Stokes operator, we may infer that any reference solution belongs to the the 
class \eqref{i12} on any compact subinterval of $(T^-, T^+)$, in particular, any reference solution is a strong solution of the MOB system. In addition, using a standard bootstrap argument, we conclude 
\begin{align}  
\sup_{t \in [0,T]} &\left( \| \vu (t, \cdot) \|_{W^{2 - \frac{2}{q}, q}(Q; \R^d)}  + \| \Theta (t, \cdot)  \|_{W^{2 - \frac{2}{q}, q}(Q) }   \right)\br  &\leq 
C \left( \mathcal{E},\| \vtB \|_{C^2(\Ov{Q})},\ \| G \|_{W^{1,\infty} (Q)},\mu^{-1}, \kappa^{-1}, q, (T^-)^{-1}\right)
\label{r2}
\end{align}
for any $q \in [1, \infty)$, and any $0 < T < T^+$, where the quantity $C$ is bounded for bounded arguments.

We conclude that the bound \eqref{r1} is always satisfied by any global 
weak solution if $d=2$, while boundedness of the velocity for $d=3$, though 
physically rather obvious, is one of the major problems in fluid mechanics, cf. Fefferman \cite{Feff}. In view of perspective applications 
of our results in meteorology, hypothesis \eqref{r1} does not seem exceedingly restrictive.

\subsection{Synchronized solution}

We start by introducing the interpolation operators, cf. Azouani, Olson, 
and Titi \cite{AzOlTi}. We consider a family of projections, 
\begin{equation} \label{r3}
I_\delta: L^2(Q) \to L^2(Q) 
\ \mbox{orthogonal projections}, \ I_\delta (w) \to w 
\ \mbox{as}\ \delta \to 0 \ \mbox{for any}\ 
w \in L^2(Q).
\end{equation}
We remark that in what follows we will not distinguish between the scalar form of $I_\delta$ when applied to the temperature and its vector valued variant applied component wise to the velocity field.

The deterministic synchronized MOB system is considered on the time interval $(0, T^+)$ and the nudging forces are active in the time lap $(0,T)$ with $0<T<T^+$. The overall strategy of the data assimilation method can be described as follows. The final objective is to recover 
the exact value of the {\it a priori} unknown observed solution $(\vu, \Theta)$ in the \emph{prediction period} $(T, T^+)$ knowing only 
its approximate values (measurements) $I_\delta(\vu)$, $I_\delta(\Theta)$ in the \emph{data sampling period} $(0,T)$, $T < T^+$, with typically 
$T^+ - T >> T$. 
This goal is achieved by 
solving the synchronized problem in the  time interval $(0,T^+)$, with suitable values of the nudging parameter $\Lambda$. 

Thus, we have the following deterministic syncronized MOB system:
 \begin{align}
\Div \tvu &= 0, \br
\D \tvu + \Div (\tvu \otimes \tvu) \dt + \Grad \Pi \dt &= 
\mu \Del \tvu \dt 	- \tTh \Grad G \dt - \vtB \Grad G \dt \br 
&- \Lambda \left( I_\delta (\tvu) - I_\delta (\vu) \right) \mathds{1}_{t \in [0,T]}\dt + \Lambda \vc{e}_{D, \vu} \mathds{1}_{t \in [0,T]}  \dt
+ \Lambda \mathds{1}_{t \in [0,T]} \D W_{\vu},\br 
\D M(\tTh) + \Div (\tvu \tTh) \dt  &= \kappa \Del \tTh \dt - \Div (\tvu (\vtB + aG)) \dt \br
&- \Lambda \left( I_\delta (\tTh) - I_\delta (\Theta) \right) \mathds{1}_{t \in [0,T]}\dt + \Lambda{e}_{D,\theta} \mathds{1}_{t \in [0,T]}  \dt
+ \Lambda \mathds{1}_{t \in [0,T]} \D W_{\Theta},
\label{r6}	
\end{align}
with the boundary conditions 
\begin{equation} \label{r7}
	\tvu|_{\partial Q} = 0, \ \tTh|_{\partial Q} = 0.
\end{equation}
Here $\Lambda$ is a relaxation parameter (nudging) that will be chosen later, forcing the syncronized solution to be close to the observed solution $(\vu, \Theta)$.

As in \cite{BeOlTi}, in what follows we consider that the actual interpolated measurements of the solution $(\vu, \Theta)$ contain random errors, this instead of $(I_\delta(\vu), I_\delta(\Theta))$ we have:
$$
(\bar I_\delta(\vu), \bar I_\delta(\Theta))=(I_\delta(\vu), I_\delta(\Theta))+(\xi_1(t), \xi_2(t))
$$
where the error $(\xi_1, \xi_2)$ is decomposed into a deterministic part $(e_{D, \vu}, e_{D, \Theta}) $ and a random part expressed in terms of cylindrical Wiener processes. On the deterministic component of the observation error we suppose that they are bounded possibly random variables.

Thus, we introduce two cylindrical Wiener processes, 
\begin{equation} \label{r4}
W_\vu = \vc{v} \cdot \chi_{\vu} = 
\sum_{k=1}^\infty v_k \chi_{\vu,k} ,\ W_\Theta = \vc{h} \cdot \chi_{\Theta, k} = \sum_{k=1}^\infty h_k \chi_{\Theta,k}, 
\end{equation}
where $( \chi_{\vu,k} )_{k=1}^\infty$, 
$(\chi_{\Theta,k} )_{k=1}^\infty$ are mutually independent Wiener processes, with variance $\sigma_\vu$, $\sigma_\Theta$, respectively. 
The	diffusion coefficients are deterministic functions of the spatial 
variable $x$, 
\begin{equation} \label{r5}
v_k \in L^\infty(Q; \R^d),\ 
h_k \in L^\infty(Q), \ \sum_{k=1}^\infty 
\left( \| v_k \|_{L^\infty(Q; \R^d)}^2 + \| h_k \|_{ L^\infty(Q)}^2 \right) \leq E_S. 
\end{equation}

The stochastic \emph{synchronized MOB system} reads as follows:

\begin{align}
\Div \tvu &= 0, \br
\D \tvu + \Div (\tvu \otimes \tvu) \dt + \Grad \Pi \dt &= 
\mu \Del \tvu \dt 	- \tTh \Grad G \dt - \vtB \Grad G \dt \br 
&- \Lambda \left( I_\delta (\tvu) - I_\delta (\vu) \right) \mathds{1}_{t \in [0,T]}\dt + \Lambda \vc{e}_{D, \vu} \mathds{1}_{t \in [0,T]}  \dt
+ \Lambda \mathds{1}_{t \in [0,T]} \D W_{\vu},\br 
\D M(\tTh) + \Div (\tvu \tTh) \dt  &= \kappa \Del \tTh \dt - \Div (\tvu (\vtB + aG)) \dt \br
&- \Lambda \left( I_\delta (\tTh) - I_\delta (\Theta) \right) \mathds{1}_{t \in [0,T]}\dt + \Lambda{e}_{D,\theta} \mathds{1}_{t \in [0,T]}  \dt
+ \Lambda \mathds{1}_{t \in [0,T]} \D W_{\Theta},
\label{r6}	
\end{align}
with the boundary conditions 
\begin{equation} \label{r7}
	\tvu|_{\partial Q} = 0, \ \tTh|_{\partial Q} = 0.
\end{equation}
The synchronized system will be solved in the time interval $[0, T^+]$ supplemented with arbitrary deterministic initial conditions
\begin{equation} \label{r8}
\tvu(0, \cdot) = \tvu_0,\ \tTh(0, \cdot) = \tTh_0.
\end{equation}

%The random observation 
%error is represented by the stochastic forcing imposed in the synchronized system \eqref{r6}. 

\subsection{Weak martingale solutions of the synchronized system}

The existence of (stochastically) strong solutions to the synchronized system for $d=3$ is not known, while
the case $d=2$ could possibly be handled, cf. e.g. Flandoli and Maslowski \cite{FlaMas}.
To remedy this obstacle, 
we introduce the concept of \emph{weak martingale solution}. 

\begin{Definition}[{\bf Weak martingale solution}] \label{rD1}
	
The quantity $\left[ \left(\Omega, \mathfrak{F}, (\mathfrak{F}_t)_{t \geq 0}, \prst \right), \tvu, \tTh, 
W_\vu, W_\Theta \right]$ is called 
\emph{weak martingale solution} of the 
synchronized MOB system \eqref{r6}--\eqref{r8}, if the following holds:
\begin{itemize}
	\item $\left(\Omega, \mathfrak{F}, (\mathfrak{F}_t)_{t \geq 0}, \prst \right)$ is a stochastic basis with a complete right continuous filtration.  
	\item
	\[ 
	W_\vu = \vc{v} \cdot \chi_{\vu} = 
	\sum_{k=1}^\infty v_k \chi_{\vu,k} ,\ W_\Theta = \vc{h} \cdot \chi_{\Theta, k} = \sum_{k=1}^\infty h_k \chi_{\Theta,k}, 
	\]
where $(\chi_{\vu,k})_{k=1}^\infty$, $(\chi_{\Theta,k})_{k=1}^\infty$ are $(\mathfrak{F}_t)$ adapted cylindrical Wiener processes, with 
variance $\sigma_\vu$, $\sigma_\Theta$, respectively.
\item 
The velocity $\tvu$ and the temperature $\tTh$ are $(\mathfrak{F}_t)$ 
progressively measurable random processes weakly-continuous in $[0, T^+]$  ranging in the 
space $W^{-\ell,2}(Q)$, $\ell > 2$.
\item 
\begin{equation} \label{r9}
\intQ{ \tvu (t, \cdot) \cdot \Grad \phi } = 0 
\ \mbox{for all}\ t \in [0,T^+ ) \ \mbox{and any} 
\ \phi \in C^1(\Ov{Q})	, \prst \text{-}\rm{a.s.}
\end{equation}	
\item 
The integral identity 
\begin{align}
- \int_0^{\tau} &\intQ{ \tvu \cdot \partial_t \bfphi } \dt +
\intQ{ \tvu \cdot \bfphi } \Big]_0^\tau= 
\int_0^{\tau} \intQ{ \tvu \otimes \tvu : \Grad \bfphi } \dt \br&- 
\mu \int_0^{\red{\tau}} \intQ{ \Grad \tvu : \Grad \bfphi } \dt  
- \int_0^{\red{\tau}} \intQ{ (\tTh + \vtB) \Grad G \cdot \bfphi } \dt \br
&- \Lambda \int_0^{\tau} \intQ{ \Big( I_\delta (\tvu) - I_\delta (\vu) \Big) \cdot \bfphi } \mathds{1}_{t \in [0,T]}\dt \br 
&+ \Lambda \int_0^{\tau} \intQ{ \vc{e}_{D,\vu} \cdot  \bfphi }  \mathds{1}_{t \in [0,T]}  \dt
+ \Lambda \sum_{k=1}^\infty \int_0^{\red{\tau}} \left( \intQ{ \bfphi \cdot 
\vc{v}_k } \right) \mathds{1}_{t \in [0,T]} \D \chi_{\vu,k}
\label{r10}
\end{align}
holds $\prst$ -a.s. for any $0 \leq \tau \leq T^+$ and any test function $\bfphi \in C^1_c([0,T^+) \times Q); \R^d$,\ 
$\Div \bfphi = 0$.

\item The integral identity
\begin{align}
	- \int_0^{\tau} &\intQ{ M(\tTh) \partial_t \varphi } \dt - 
	\intQ{ M(\tTh_0)  \varphi(0,\cdot) } = 
	\int_0^{\tau} \intQ{ \tTh \tvu \cdot \Grad \varphi } \dt \br&- 
	\kappa \int_0^{\tau} \intQ{ \Grad \tTh \cdot\Grad \varphi } \dt  
	+ \int_0^{\tau} \intQ{  (\tvu (\vtB + aG))   \cdot \Grad \varphi } \dt \br
	&- \Lambda \int_0^{\tau} \intQ{ \Big( I_\delta (\tTh) - I_\delta (\Theta) \Big)  \varphi } \mathds{1}_{t \in [0,T]}\dt \br 
	&+ \Lambda \int_0^{\tau} \intQ{ {e}_{D,\Theta} \cdot  \varphi }  \mathds{1}_{t \in [0,T]}  \dt
	+ \Lambda \sum_{k=1}^\infty \int_0^{\tau} \left( \intQ{ \varphi  
		h_k } \right) \mathds{1}_{t \in [0,T]} \D \chi_{\Theta,k}
	\label{r11}
\end{align}
holds  $\prst$ -a.s. for any $0 \leq \tau \leq T^+$ and any test function $\varphi \in C^1_c([0, T^+) \times Q)$.

\item The kinetic energy inequality 
\begin{align}
- \int_0^{T^+} &\partial_t \psi \intQ{ \frac{1}{2} |\tvu|^2 } \dt
- \psi(0) \intQ{ \frac{1}{2} |\tvu_0|^2 } + \mu 
\int_0^{T^+} \psi \intQ{ |\Grad \tvu|^2 } \dt \br
&+ \Lambda \int_0^{T^+} \psi \intQ{ \Big( I_\delta(\tvu) - I_\delta (\vu) \Big)\mathds{1}_{t \in [0,T]}  \cdot \tvu } \dt \br 
&\leq - \int_0^{T^+} \psi \intQ{ (\tTh + \vtB) \Grad G \cdot \tvu } \dt
+ \Lambda \int_0^{T^+} \psi \intQ{ \vc{e}_{D, \vu} \cdot \tvu }\mathds{1}_{t \in [0,T]} \dt \br 
& \quad + \frac{1}{2} \Lambda^2 \sigma^2_\vu \sum_{k=1}^\infty \int_0^{T^+} \psi \left( \intQ{  
	|\vc{v}_k|^2 } \right) \mathds{1}_{t \in [0,T]} \dt              \br
& \quad + \Lambda \sum_{k=1}^\infty \int_0^{T^+} \psi \left( \intQ{ \tvu \cdot 
	\vc{v}_k } \right) \mathds{1}_{t \in [0,T]} \D \chi_{\vu,k}
\label{r12}	
\end{align}
holds $\prst$-a.s. for any $\psi \in C^1_c[0, T^+)$, $\psi \geq 0$.
\item 
The thermal energy balance 
\begin{align}
	- \int_0^{T^+} &\partial_t \psi \intQ{ \frac{1}{2} M(\tTh) \tTh } \dt
	- \psi(0) \intQ{ \frac{1}{2} M(\tTh_0) \tTh_0 } + \kappa 
	\int_0^{T^+} \psi \intQ{ |\Grad \tTh|^2 } \dt \br
	&+ \Lambda \int_0^{T^+} \psi \intQ{ \Big( I_\delta(\tTh) - I_\delta (\Theta) \Big) \mathds{1}_{t \in [0,T]} \cdot \tTh } \dt \br 
	&= \int_0^{T^+} \psi \intQ{ (\vtB + a G) \tvu \cdot \Grad \tTh } \dt
	+ \Lambda \int_0^{T^+} \psi \intQ{ {e}_{D, \Theta} \tTh }\mathds{1}_{t \in [0,T]} \dt \br 
	& \quad + \frac{1}{2} \Lambda^2 \sigma^2_\Theta \sum_{k=1}^\infty \int_0^{T^+} \psi \left( \intQ{  
		h_k M^{-1} (h_k) } \right) \mathds{1}_{t \in [0,T]} \dt              \br
	& \quad + \Lambda \sum_{k=1}^\infty \int_0^{T^+} \psi \left( \intQ{ \tTh h_k } \right) \mathds{1}_{t \in [0,T]} \D \chi_{\vu,k}
	\label{r13}	
\end{align}
holds $\prst$-a.s. for any $\psi \in C^1_c[0,T^+)$.
\end{itemize}

\end{Definition}

%The quantities $\vc{e}_{D, \vu}$, $e_{D, \Theta}$ represent a deterministic component of the observation error and we always suppose they are bounded 
%possibly random variables. 
The \emph{existence} of global in time weak--martingale solutions to the synchronized system \eqref{r6}, \eqref{r7}, with the 
initial (deterministic) data \eqref{r8} can be shown by the application of the stochastic compactness method exactly as in \cite{FeiPet}, cf. also the monograph 
\cite{BrFeHobook}.

\section{Main result}
\label{M}

Having collected the necessary preliminary material we are ready to state our main result. 

\begin{Theorem}[\bf Continuous data assimilation] \label{TM1}
Let $Q \subset R^d$, $d = 2,3$ be a bounded domain of class $C^{2, \alpha}$, $\alpha > 0$. Let $(\vu, \Theta)$ be a solution of the MOB system 
\eqref{i6}, \eqref{i7} in the time interval $(T^-, T^+)$ satisfying \eqref{r1}. Let $(\tvu, \tvT)$ be a weak martingale solution of the synchronized 
system \eqref{r6}--\eqref{r8} in $(0,T^+)$, $T^+ > T$ in the sense of Definition \ref{rD1}. 

Then for any $\gamma > 0$, there exist $\Lambda_0 > 0$, $\delta_0 > 0$ depending only on $\mu$, $\kappa$, and the specific form of the 
interpolation operators $I_\delta$ such that
\begin{align}
&\expe{ \left\| (\vu - \tvu)(\tau, \cdot) \right\|^2_{L^2(Q; \R^d) } + \left\| (\Theta - \tvT)(\tau, \cdot) \right\|^2_{L^2(Q) }} \br
&\quad \aleq \exp(K(\tau - T^+)) \left[ \exp(- \gamma T) \max \left\{ \mathcal{E}, \| \vu_0 \|_{L^2(Q; \R^3)}, \| \tvT_0 \|_{L^2(Q)} \right\} \right] \br
&\quad + \exp(K(\tau - T^+)) \Lambda^2 \int_0^T \intQ{ \expe{ \left[ |\vc{e}_{D, \vu}|^2 + |e_{D, \Theta} |^2 \right] } } \dt \br 
&\quad + \exp(K(\tau - T^+)) \Lambda^2 \left[  \frac{\sigma^2_\vu}{2} \sum_{k=1}^\infty \int_0^{T} \left( \intQ{  
	|\vc{v}_k|^2 } \right)  \dt  +  \frac{\sigma^2_\Theta}{2} \sum_{k=1}^\infty \int_0^{T} \psi \left( \intQ{  
	h_k M^{-1} (h_k) } \right)  \dt \right]
\label{M1}
\end{align}	
for all $\tau \in (T; T^+)$ and any $\Lambda \geq \Lambda_0$, $0 < \delta \leq \delta_0$, where $K$ depends only on $\mathcal{E}$, $\mu$, and $\kappa$.
\end{Theorem}	

A short inspection of formula \eqref{M1} reveals that the first term on the right--hand side can be made small by choosing 
$\gamma > 0$ large enough. This in turn forces $\Lambda_0$ to be large increasing the other two integrals that represent 
the approximation error. 

The rest of the paper is devoted to the proof of Theorem \ref{TM1}.

\section{Relative energy inequality}
\label{e}

The relative energy associated to the synchronized MOB system is the quantity 
\[
E \left( \tvu, \tTh \Big| \vc{w}, \mathcal{T} \right) = \intQ{
\frac{1}{2} \left( |\tvu - \vc{w}|^2 + (M(\tTh) - M(\mathcal{T}))(\tTh - \mathcal{T} ) \right)}.
\]
Seeing that 
\[
E \left( \tvu, \tTh \Big| \vc{w}, \mathcal{T} \right) = \intQ{
\frac{1}{2} |\tvu|^2 + \frac{1}{2} M(\tTh) \tTh  
- \tvu \cdot \vc{w} - \frac{1}{2}( M(\tTh) \mathcal{T} + 
M(\mathcal{T}) \tTh )  + \frac{1}{2} |\vc{w}|^2 + \frac{1}{2}M (\mathcal{T}) \mathcal{T}},
\]
we may trace the time evolution of the relative energy using the 
weak formulation \eqref{r10} -- \eqref{r13} as soon as 
$\vc{w}$ and $\mathcal{T}$ are sufficiently smooth deterministic functions satisfying the compatibility conditions 
\begin{equation} \label{e1}
\Div \vc{w} = 0,\ \vc{w}|_{\partial Q} = 0,\ \mathcal{T}|_{\partial Q} = 0.
\end{equation}	

Indeed, plugging $\vc{w}$ as a test function in the momentum balance 
\eqref{r10} we obtain 
\begin{align} 
	\left[ \intQ{ \tvu \cdot \vc{w} } \right]_{t=0}^{t = \tau} 
	 &=  \int_0^{\tau} \intQ{ \tvu \cdot \partial_t \vc{w} } \dt +
	\int_0^{\tau} \intQ{ \tvu \otimes \tvu : \Grad \vc{w} } \dt \br&- 
	\mu \int_0^{\tau} \intQ{ \Grad \tvu : \Grad \vc{w} } \dt  
	- \int_0^{\tau} \intQ{ (\tTh + \vtB) \Grad G \cdot \vc{w}  } \dt \br
	&- \Lambda \int_0^{\tau} \intQ{ \Big( I_\delta (\tvu) - I_\delta (\vu) \Big) \cdot \vc{w} } \mathds{1}_{t \in [0,T]}\dt \br 
	&+ \Lambda \int_0^{\tau} \intQ{ \vc{e}_{D,\vu} \cdot  \vc{w} }  \mathds{1}_{t \in [0,T]}  \dt
	+ \Lambda \sum_{k=1}^\infty \int_0^{\tau} \left( \intQ{ \vc{w} \cdot 
		\vc{v}_k } \right) \mathds{1}_{t \in [0,T]} \D \chi_{\vu,k}
	\label{e2}
\end{align}
for any $0 \leq \tau < T^+$ $\prst$-a.s.

Repeating the same argument with the heat equation \eqref{r11} we obtain 
\begin{align}
	  \left[ 
	\intQ{ M(\tTh) \mathcal{T} } \right]_{t = 0}^{t = \tau} &= 
	\int_0^{\tau} \intQ{ M(\tTh) \partial_t \mathcal{T} } \dt +
	\int_0^{\tau} \intQ{ \tTh \tvu \cdot \Grad \mathcal{T} } \dt \br&- 
	\kappa \int_0^{\tau} \intQ{ \Grad \tTh \cdot\Grad \mathcal{T} } \dt  
	+ \int_0^{\tau} \intQ{  (\tvu (\vtB + aG))   \cdot \Grad \mathcal{T} } \dt \br
	&- \Lambda \int_0^{\tau} \intQ{ \Big( I_\delta (\tTh) - I_\delta (\Theta) \Big)  \mathcal{T} } \mathds{1}_{t \in [0,T]}\dt \br 
	&+ \Lambda \int_0^{\tau} \intQ{ {e}_{D,\Theta} \cdot  \mathcal{T} }  \mathds{1}_{t \in [0,T]}  \dt
	+ \Lambda \sum_{k=1}^\infty \int_0^{\tau} \left( \intQ{ \mathcal{T}  
		h_k } \right) \mathds{1}_{t \in [0,T]} \D \chi_{\Theta,k}
	\label{e3}
\end{align}
for any $0 \leq \tau < T^+$ $\prst$-a.s.

Replacing $\vc{w}$ by $\psi \vc{w}$, $\psi \in C^1_c [0,T^+)$, $\psi \geq 0$ in \eqref{e2}, we substract the resulting equation from  the kinetic energy inequality \eqref{r12} and obtain the \emph{relative energy inequality} 
for the kinetic energy:
\begin{align}
	- \int_0^{T^+} &\partial_t \psi \intQ{ \frac{1}{2} \left| \tvu - \vc{w}  \right|^2 } \dt
 + \mu 
	\int_0^{T^+} \psi \intQ{ \Grad \tvu \cdot (\Grad \tvu - \Grad \vc{w} ) } \dt \br
	&+ \Lambda \int_0^{T^+} \psi \intQ{ \Big( I_\delta(\tvu) - I_\delta (\vu) \Big) \cdot (\tvu - \vc{w} ) }\mathds{1}_{t \in [0,T]} \dt \br 
	&\leq 	\psi(0) \intQ{ \frac{1}{2} \left| \tvu_0 - \vc{w}(0, \cdot) \right|^2 }  
 \br
	&+ \int_0^{T^+} \psi \intQ{ (\vc{w} - \tvu ) \cdot \partial_t \vc{w} } \dt -
	\int_0^{T^+} \psi \intQ{ \tvu \otimes \tvu : \Grad \vc{w} } \dt \br  &+ \int_0^{T^+} \psi \intQ{ (\tTh + \vtB) \Grad G \cdot (\vc{w} - \tvu)  } \dt
	\br
	&
	+ \Lambda \int_0^{T^+} \psi \intQ{ \vc{e}_{D, \vu} \cdot (\tvu - \vc{w}) }\mathds{1}_{t \in [0,T]} \dt \br 
	& + \frac{1}{2} \Lambda^2 \sigma^2_\vu \sum_{k=1}^\infty \int_0^{T^+} \psi \left( \intQ{  
		|\vc{v}_k|^2 } \right) \mathds{1}_{t \in [0,T]} \dt              \br
	& + \Lambda \sum_{k=1}^\infty \int_0^{T^+} \psi \left( \intQ{ (\tvu - \vc{w}) \cdot 
		\vc{v}_k } \right) \mathds{1}_{t \in [0,T]} \D \chi_{\vu,k}
	\label{e4}	
\end{align}
for any $\psi \in C^1_c[0,T^+)$, and any ``test function'' $\vc{w}$ satisfying \eqref{e1} $\prst$-a.s.

The thermal energy balance can be handled in a similar manner. First observe that 
\[
\intQ{ M(\tTh) \mathcal{T} } = 
\intQ{ \tTh M (\mathcal{T}) }. 
\]
Consequently, replacing $\mathcal{T}$ by $\psi \mathcal{T}$ in \eqref{e4}, with $\psi \in \mathcal{C}^1_c([0, T^+))$ and combining to the thermal energy balance \eqref{r13}, we obtain the \emph{relative energy inequality} for the thermal equation:
\begin{align}
	- \int_0^{T^+} &\partial_t \psi \intQ{ \frac{1}{2} \left( M(\tTh - 
	\mathcal{T} )(\tTh - \mathcal{T}) \right) } \dt
	 + \kappa 
	\int_0^{T^+} \psi \intQ{ \Grad \tTh \cdot (\Grad \tTh - \Grad \mathcal{T} ) } \dt \br
	&+ \Lambda \int_0^{T^+} \psi \intQ{ \Big( I_\delta(\tTh) - I_\delta (\Theta) \Big) ( \tTh - \mathcal{T})  } \dt \br 
	&= \psi(0) \intQ{ \frac{1}{2} \left( ( M(\tTh_0 - \mathcal{T}(0, \cdot))( \tTh_0 - \mathcal{T}(0, \cdot)) \right) } \br
	&-\int_0^{T^+} \psi \intQ{ ( M(\tTh) - M(\mathcal{T}) ) \partial_t \mathcal{T} } \dt -
	\int_0^{T^+} \psi \intQ{ \tTh \tvu \cdot \Grad \mathcal{T} } \dt	\br
	&+\int_0^{T^+} \psi \intQ{ (\vtB + a G) \tvu \cdot (\Grad \tTh - \Grad \mathcal{T}) } \dt
	+ \Lambda \int_0^{T^+} \psi \intQ{ {e}_{D, \Theta} (\tTh - \mathcal{T}) }\mathds{1}_{t \in [0,T]} \dt \br 
	& + \frac{1}{2} \Lambda^2 \sigma^2_\Theta \sum_{k=1}^\infty \int_0^{T^+} \psi \left( \intQ{  
		h_k M^{-1} (h_k) } \right) \mathds{1}_{t \in [0,T]} \dt              \br
	&+ \Lambda \sum_{k=1}^\infty \int_0^{T^+} \psi \left( \intQ{ (\tTh - \mathcal{T}) h_k } \right) \mathds{1}_{t \in [0,T]} \D \chi_{\red{\theta},k}
	\label{e5}	
\end{align}
It is worth noting that the relations \eqref{e4}, \eqref{e5} hold for 
any pair of differentiable deterministic functions $(\vc{w}, \mathcal{T})
$ satisfying the compatibility conditions \eqref{e1}, in particular for any solution $(\vu, \Theta)$ of the observed system. 

\section{Approximating the reference solution}
\label{a} 

As an application of the relative energy inequalities \eqref{e4}, \eqref{e5}, we obtain a bound on the distance of the observed and synchronized solutions. As the reference solution is smooth, it can be used as test functions $\vc{w} = \vu$, $\mathcal{T} = \Theta$ in 
\eqref{e3}, \eqref{e4}:
\begin{align}
	- \int_0^{T^+} &\partial_t \psi \intQ{ \frac{1}{2} \left| \tvu - \vu  \right|^2 } \dt
	+ \mu 
	\int_0^{T^+} \psi \intQ{ \Grad \tvu \cdot (\Grad \tvu - \Grad \vu ) } \dt \br
	&+ \Lambda \int_0^{T^+} \psi \intQ{ \Big| I_\delta(\tvu) - I_\delta (\vu) \Big|^2 }\mathds{1}_{t \in [0,T]} \dt \br 
	&\leq 	\psi(0) \intQ{ \frac{1}{2} \left| \tvu_0 - \vu(0, \cdot) \right|^2 }  
	\br
	&+ \int_0^{T^+} \psi \intQ{ (\vu - \tvu ) \cdot \partial_t \vu } \dt -
	\int_0^{T^+} \psi \intQ{ \tvu \otimes \tvu : \Grad \vu } \dt \br  &+ \int_0^{T^+} \psi \intQ{ (\tTh + \vtB) \Grad G \cdot (\vu - \tvu)  } \dt
	\br
	&
	+ \Lambda \int_0^{T^+} \psi \intQ{ \vc{e}_{D, \vu} \cdot (\tvu - \vu) }\mathds{1}_{t \in [0,T]} \dt \br 
	& + \frac{1}{2} \Lambda^2 \sigma^2_\vu \sum_{k=1}^\infty \int_0^{T^+} \psi \left( \intQ{  
		|\vc{v}_k|^2 } \right) \mathds{1}_{t \in [0,T]} \dt              \br
	& + \Lambda \sum_{k=1}^\infty \int_0^{T^+} \psi \left( \intQ{ (\tvu - \vu) \cdot 
		\vc{v}_k } \right) \mathds{1}_{t \in [0,T]} \D \chi_{\vu,k},
	\label{a1}	
\end{align}
and 
\begin{align}
	- \int_0^{T^+} &\partial_t \psi \intQ{ \frac{1}{2} \left( M(\tTh - 
		\Theta )(\tTh - \Theta) \right) } \dt
	+ \kappa 
	\int_0^{T^+} \psi \intQ{ \Grad \tTh \cdot (\Grad \tTh - \Grad \Theta ) } \dt \br
	&+ \Lambda \int_0^{T^+} \psi \intQ{ \Big| I_\delta(\tTh) - I_\delta (\Theta) \Big|^2  } \dt \br 
	&= \psi(0) \intQ{ \frac{1}{2} \left( ( M(\tTh_0 - \Theta(0, \cdot))( \tTh_0 - \Theta(0, \cdot)) \right) } \br
	&-\int_0^{T^+} \psi \intQ{ ( M(\tTh) - M(\Theta) ) \partial_t \Theta } \dt -
	\int_0^{T^+} \psi \intQ{ \tTh \tvu \cdot \Grad \Theta } \dt	\br
	&+\int_0^{T^+} \psi \intQ{ (\vtB + a G) \tvu \cdot (\Grad \tTh - \Grad \Theta) } \dt
	+ \Lambda \int_0^{T^+} \psi \intQ{ {e}_{D, \Theta} (\tTh - \Theta) }\mathds{1}_{t \in [0,T]} \dt \br 
	& + \frac{1}{2} \Lambda^2 \sigma^2_\Theta \sum_{k=1}^\infty \int_0^{T^+} \psi \left( \intQ{  
		h_k M^{-1} (h_k) } \right) \mathds{1}_{t \in [0,T]} \dt              \br
	&+ \Lambda \sum_{k=1}^\infty \int_0^{T^+} \psi \left( \intQ{ (\tTh - \Theta) h_k } \right) \mathds{1}_{t \in [0,T]} \D \chi_{\vu,k}.
	\label{a2}	
\end{align}

Next, we exploit the fact that $\vu$, $\Theta$ satisfy the field equations in the MOB system \eqref{i6} to rewrite \eqref{a1} in the form 
\begin{align}
	- \int_0^{T^+} &\partial_t \psi \intQ{ \frac{1}{2} \left| \tvu - \vu  \right|^2 } \dt \br
	&+ \mu 
	\int_0^{T^+} \psi \intQ{ |\Grad \tvu - \Grad \vu |^2 } \dt + \Lambda \int_0^{T^+} \psi \intQ{ \Big| I_\delta(\tvu) - I_\delta (\vu) \Big|^2 }\mathds{1}_{t \in [0,T]} \dt \br 
	&\leq 	\psi(0) \intQ{ \frac{1}{2} \left| \tvu_0 - \vu(0, \cdot) \right|^2 }  
	\br
	&+\int_0^{T^+} \psi \intQ{ (\vu - \tvu) \otimes (\tvu - \vu) : \Grad \vu } \dt + \int_0^{T^+} \psi \intQ{ (\tTh - \Theta) \Grad G \cdot (\vu - \tvu)  } \dt
	\br
	&
	+ \Lambda \int_0^{T^+} \psi \intQ{ \vc{e}_{D, \vu} \cdot (\tvu - \vu) }\mathds{1}_{t \in [0,T]} \dt \br 
	& + \frac{1}{2} \Lambda^2 \sigma^2_\vu \sum_{k=1}^\infty \int_0^{T^+} \psi \left( \intQ{  
		|\vc{v}_k|^2 } \right) \mathds{1}_{t \in [0,T]} \dt              \br
	& + \Lambda \sum_{k=1}^\infty \int_0^{T^+} \psi \left( \intQ{ (\tvu - \vu) \cdot 
		\vc{v}_k } \right) \mathds{1}_{t \in [0,T]} \D \chi_{\vu,k},
	\label{a3}	
\end{align}
where we have used 
\[
\intQ{ \vu \otimes \tvu : \Grad \vu } = 
\intQ{ \frac{1}{2}\tvu \cdot \Grad |\vu|^2 }   = 0.
\]

Applying the same treatment to the thermal energy inequality \eqref{a2} 
we obtain
\begin{align}
	- \int_0^{T^+} &\partial_t \psi \intQ{ \frac{1}{2} \left( M(\tTh - 
		\Theta )(\tTh - \Theta) \right) } \dt
	+ \kappa 
	\int_0^{T^+} \psi \intQ{  | \Grad \tTh - \Grad \Theta |^2 } \dt \br
	&+ \Lambda \int_0^{T^+} \psi \intQ{ \Big| I_\delta(\tTh) - I_\delta (\Theta) \Big|^2  }\mathds{1}_{t \in [0,T]} \dt \br 
	&= \psi(0) \intQ{ \frac{1}{2} \left( ( M(\tTh_0 - \Theta(0, \cdot))( \tTh_0 - \Theta(0, \cdot)) \right) } \br
	&+ \int_0^{T^+} \psi \intQ{ ( \tvu -  \vu) \Theta \cdot \Grad (\tTh - \Theta)  \Big) } \dt 	\br
	&+\int_0^{T^+} \psi \intQ{ (\vtB + a G) (\tvu - \vu) \cdot (\Grad \tTh - \Grad \Theta) } \dt \br
	&+ \Lambda \int_0^{T^+} \psi \intQ{ {e}_{D, \Theta} (\tTh - \Theta) }\mathds{1}_{t \in [0,T]} \dt \br 
	& + \frac{1}{2} \Lambda^2 \sigma^2_\Theta \sum_{k=1}^\infty \int_0^{T^+} \psi \left( \intQ{  
		h_k M^{-1} (h_k) } \right) \mathds{1}_{t \in [0,T]} \dt              \br
	&+ \Lambda \sum_{k=1}^\infty \int_0^{T^+} \psi \left( \intQ{ (\tTh - \Theta) h_k } \right) \mathds{1}_{t \in [0,T]} \D \chi_{\vu,k}.
	\label{a4}	
\end{align}

Going back to \eqref{a3} we can control 
\begin{equation} \label{a5}
\left| \intQ{ (\vu - \tvu)\cdot \Grad \vu \cdot (\tvu - \vu) } \right| \leq 
\| \vu \|_{L^\infty(Q; \R^d)} \| \Grad (\vu -\tvu) \|_{L^2(Q; \R^{d \times d})} \| 
\vu - \tvu \|_{L^2(Q; \R^d)}, 
\end{equation}	
and, similarly,
\begin{equation} \label{a6}
\left| \intQ{ (\tTh - \Theta) \Grad G \cdot (\vu - \tvu) }\right| \leq 
	\| \Grad G \|_{L^\infty(Q; \R^d)} \| \tTh - \Theta \|_{L^2(Q)} \| 
	\vu - \tvu \|_{L^2(Q; \R^d)}.
\end{equation}

Next, the integrals on the right--hand side of 
\eqref{a4} are bounded as 
\begin{align}
\left| \intQ{ \Theta (\tvu - \vu) \cdot \Grad (\tTh - \Theta) } \right|
+ \left|\intQ{ (\vtB + a G) (\tvu - \vu) \cdot (\Grad \tTh - \Grad \Theta) }   \right|	\br 
\leq \Big( \| \Theta \|_{L^\infty(Q)} + \| \vtB \|_{L^\infty(Q)} + 
|a| \| G \|_{L^\infty(Q)} \Big) \| \tvu - \vu \|_{L^2(Q; \R^d)} 
\| \Grad (\tTh - \Theta) \|_{L^2(Q; \R^d)} . 
\label{a7}	
\end{align}

In view of the approximation property \eqref{r3} of the interpolation operators, it is a routine matter to show the following result: 

\begin{Lemma} \label{aL1}
Suppose the interpolation operators $I_\delta$ satisfy \eqref{r3}. 

Then for any $\beta > 0$, there exists $\Lambda_0 > 0$, $\delta_0 > 0$ such that 
\begin{align}
\mu \| \Grad \vc{w} \|_{L^2(Q; \R^{d \times d})}^2 + 
\Lambda \| I_\delta (\vc{w}) \|^2_{L^2(Q; \R^d)} &\geq \frac{\mu}{2} \|\Grad \vc{w} \|_{L^2(Q; \R^{d \times d})}^2 + \beta \| \vc{w} \|^2_{L^2(Q; \R^d)}, \label{a8}\\ 
\kappa \| \Grad \mathcal{T} \|_{L^2(Q; \R^{d})}^2 + 
\Lambda \| I_\delta (\mathcal{T}) \|^2_{L^2(Q)} &\geq \frac{\kappa}{2} \|\Grad \mathcal{T} \|_{L^2(Q; \R^{})}^2 + \beta \| \mathcal{T} \|^2_{L^2(Q)} \label{a9}
\end{align}
for all $\Lambda \geq \Lambda_0$, $0 < \delta < \delta_0$. The values 
of $\Lambda_0$, $\delta_0$ depend on $\beta$, $\mu$, and $\kappa$.
\end{Lemma}	

Summing up the inequalities \eqref{a3}, \eqref{a4}, and using the estimates  
\eqref{a5}--\eqref{a7} together with Lemma \ref{aL1}, we conclude that 
for any $\gamma > 0$, there exists $\Lambda_0 > 0$, $\delta_0 > 0$ such that 
\begin{align}
- & \int_0^{T} \partial_t \psi \intQ{ \frac{1}{2}\left[ | \tvu - \vu |^2 + M(\tTh - \Theta)(\tTh - \Theta))  \right] } \dt \br 	
& + \gamma \int_0^{T} \psi \intQ{\left[ | \tvu - \vu |^2 + M(\tTh - \Theta)(\tTh - \Theta))  \right] } \dt \br 
&\leq \psi(0) \intQ{\frac{1}{2}\left[ | \tvu_0 - \vu(0, \cdot) |^2 + M(\tTh_0 - \Theta(0, \cdot))(\tTh_0 - \Theta(0, \cdot))  \right]  } \br
& \Lambda^2 \int_0^{T} \psi \intQ{ \left[ |\vc{e}_{D, \vu}|^2 + |e_{D, \Theta} |^2 \right] } \dt \br 
 & + \frac{1}{2} \Lambda^2 \sigma^2_\vu \sum_{k=1}^\infty \int_0^{T} \psi \left( \intQ{  
 	|\vc{v}_k|^2 } \right)  \dt  +  \frac{1}{2} \Lambda^2 \sigma^2_\Theta \sum_{k=1}^\infty \int_0^{T} \psi \left( \intQ{  
 	h_k M^{-1} (h_k) } \right)  \dt              \br          \br
 & + \Lambda \sum_{k=1}^\infty \int_0^{T} \psi \left( \intQ{ (\tvu - \vu) \cdot 
 	\vc{v}_k } \right)  \D \chi_{\vu,k} + \Lambda \sum_{k=1}^\infty \int_0^{T} \psi \left( \intQ{ (\tTh - \Theta) h_k } \right)  \D \chi_{\red{\theta},k}. 	
\label{a10}
\end{align}
for all $\psi \in C^1_c[0,T)$, $\psi \geq 0$ whenever $\Lambda \geq \Lambda_0$, $0 < \delta \leq \delta_0$ $\prst$-a.s. 
Passing to expectations, we obtain the desired conclusion
\begin{align}
- & \int_0^{T} \partial_t \psi \expe{ \intQ{ \frac{1}{2}\left[ | \tvu - \vu |^2 + M(\tTh - \Theta)(\tTh - \Theta))  \right] } } \dt \br 	
& + \gamma \int_0^{T} \psi \expe{\intQ{\left[ | \tvu - \vu |^2 + M(\tTh - \Theta)(\tTh - \Theta))  \right] }} \dt \br 
&\leq \psi(0) \expe{ \intQ{\frac{1}{2}\left[ | \tvu_0 - \vu(0, \cdot) |^2 + M(\tTh_0 - \Theta(0, \cdot))(\tTh_0 - \Theta(0, \cdot))  \right]  }} \br
&+ \Lambda^2 \int_0^{T} \psi \intQ{ \expe{ \left[ |\vc{e}_{D, \vu}|^2 + |e_{D, \Theta} |^2 \right] } } \dt
\br
 & + \frac{1}{2} \Lambda^2 \sigma^2_\vu \sum_{k=1}^\infty \int_0^{T} \psi \left( \intQ{  
 	|\vc{v}_k|^2 } \right)  \dt  +  \frac{1}{2} \Lambda^2 \sigma^2_\Theta \sum_{k=1}^\infty \int_0^{T} \psi \left( \intQ{  
 	h_k M^{-1} (h_k) } \right)  \dt
\label{a11}
\end{align}
for any $\psi \in C^1_c [0,T)$, $\psi \geq 0$.

Relation \eqref{a11} yields exponential decay of the distance of the 
reference and synchronized solutions on the time interval $[0,T)$ modulo 
the observation errors.

Integrating the same system on the time interval $[T, T^+)$ we obtain
\begin{align}
	- & \int_T^{T^+} \partial_t \psi { \intQ{ \frac{1}{2}\left[ | \tvu - \vu |^2 + M(\tTh - \Theta)(\tTh - \Theta))  \right] } } \dt \br 	
	&\leq \psi(T) { \intQ{\frac{1}{2}\left[ | \tvu(T, \cdot) - \vu(T, \cdot) |^2 + M(\tTh(T, \cdot) - \Theta(T, \cdot))(\tTh(T, \cdot) - \Theta(T, \cdot))  \right]  }} \br
	&+ K \int_T^{T^+} \psi {\intQ{\left[ | \tvu - \vu |^2 + M(\tTh - \Theta)(\tTh - \Theta))  \right] }} \dt \br 
	\label{a12}
\end{align}
for any $\psi \in C^1_c [T,T^+)$, $\psi \geq 0$ $\prst$-a.s., where $K > 0$ 
depends on $\mu$, $\kappa$, and $\mathcal{E}$.

Combining \eqref{a11}, \eqref{a12} and using the fact that 
\[
\intQ{ M(\tvT - \Theta)(\tvT - \Theta) } \approx \| \tvT - \Theta \|^2_{L^2(Q)} 
\]
we obtain \eqref{M1}. We have proved Theorem \ref{TM1}.

%\bibliography{citace}
%\bibliographystyle{plain}

\def\cprime{$'$} \def\ocirc#1{\ifmmode\setbox0=\hbox{$#1$}\dimen0=\ht0
	\advance\dimen0 by1pt\rlap{\hbox to\wd0{\hss\raise\dimen0
			\hbox{\hskip.2em$\scriptscriptstyle\circ$}\hss}}#1\else {\accent"17 #1}\fi}

\end{document}